\documentclass{amsart}

\usepackage{amsmath,amssymb,amsthm}
\usepackage{mathtools}

\newcommand{\setC}{\mathbb{C}}
\newcommand{\setR}{\mathbb{R}}

\newcommand{\conj}[1]{\overline{#1}}

\newcommand{\Ball}[1]{\mathbb{B}^{#1}}

\newcommand{\Balln}{\Ball{n}}
\newcommand{\diff}{\mathop{}\!{d}}

\DeclarePairedDelimiter\abs{\lvert}{\rvert}
\DeclarePairedDelimiter\paren{\lparen}{\rparen}
\DeclarePairedDelimiterX\innerp[2]{\langle}{\rangle}{#1, #2}
\DeclarePairedDelimiterX\ccint[2]{\lbrack}{\rbrack}{#1, #2}

\newcommand{\Betaf}[2]{\operatorname{B}\paren{#1,#2}}
\newcommand{\Gamaf}[1]{\Gamma\paren{#1}}
\newcommand{\pFq}[3]{{}_{#1}{F}_{#2}\!\paren{#3}}

\newtheorem{theorem}{Theorem}

\author{Petar Melentijević}
\thanks{The author is partially supported by MPNTR grant 174017, Serbia}
\title{Norm of the Bergman projection onto the Bloch space with $\mathcal{M}$-invariant gradient norm}

\begin{document}
 \begin{abstract}
 	The operator norm of Bergman projections $P_{\alpha}$ from $L^{\infty}(\mathbb{B}^n)$ to the Bloch space was found in \cite{KalajMarkovic2014}. In the same paper the authors made a conjecture on the norms of $P_{\alpha}$ with respect to $\mathcal{M}$-invariant gradient norm. In this paper we prove their conjecture. 
 \end{abstract} 
\maketitle 
    
 \section{Introduction}  
 \subsection{Bergman projection }
 Let $\mathbb{B}^n$ denote the unit ball in $\mathbb{C}^n$, $n \geq 1$ and let be $dv_{\alpha}$  the measure given by 
 $$d v_{\alpha} (z) = c_{\alpha} \left(1- |z|^2 \right)^{\alpha}d v(z),$$ 
 where $d v(z)$ is the Lebesgue measure on $\mathbb{B}^n$ and 
 \begin{equation}
 \label{calfa}
 c_{\alpha}= \frac{\Gamma(n+\alpha+1)}{\Gamma(\alpha+1)\pi^n}
 \end{equation} 
 is a normalizing constant i.e. $v_{\alpha}(\mathbb{B}^n)=1$. We will also use the symbol $v_n$ for the $n-$dimensional 
 Lebesgue measure at places where dimensions must be distinguished.
 
 For $\alpha>-1$ the Bergman projection operator is given by 
 $$  P_{\alpha} f(z)= c_{\alpha}\int_{\mathbb{B}^n} K_{\alpha}(z,w) f(w) dv_{\alpha}(w) , \quad  f \in L^p({\mathbb{B}^n}),\quad  1 <p \leq \infty $$
  where
   $$ K_{\alpha}(z,w)=\frac{1}{(1- \langle z,w \rangle)^{n+\alpha+1}}, \qquad z,w \in \mathbb{B}^n. $$
 
 Here $\langle z,w \rangle$ stands for scalar product given by $z_1\overline{w}_1+z_2\overline{w}_2+\dots+z_n\overline{w}_n.$ 
 These projections are among the most important operators in theory of analytic function spaces. In \cite{ForelliRudin1974}, Forelli and Rudin proved that $P_{\alpha}$ is bounded as operator from $L^p(\mathbb{B}^n)$ to Bergman space of all $p-$integrable analytic functions on $\mathbb{B}^n$ if and only if $\alpha>\frac{1}{p}-1.$They also found the exact operator norm in cases $p=1$ and $p=2.$
 Mateljević and Pavlović extended these result for $0<p<1,$ see \cite{MateljevicPavlovic1993}.
 
 The problem of finding the exact value of the operator norm of $P_{\alpha}$ on $L^p$ spaces  turned out to be quite difficult, even for $P=P_0$. In \cite{Zhu2006}, Zhu obtained asymptotically sharp two sided norm estimates , while Dostanić in \cite{M.R.Dostanic2008} gave the following estimate: $$ \frac{1}{2} \csc\frac{\pi}{p} \leq ||P||_{p} \leq \pi \csc\frac{\pi}{p}, $$
 for $1<p<\infty.$ Liu improved these estimates in \cite{Liu2015}, for the unit ball $\mathbb{B}^n.$ Also, papers \cite{Liu2017} and \cite{LiuPeralaZhou} give the estimates for the Bergman projection in the Siegel Upper-Half space and for the weighted Bergman projections in the unit disk.
 In recent years, there has been increasing interest in studying projections of this type in various spaces. See also \cite{Kaptanoglu2005}, \cite{Vujadinovic2013}, \cite{Markovic2014}.
 
 Here, we investigate the operator norm of $P_{\alpha}:L^{\infty} \to  \mathcal{B}$ with $\mathcal{M}$-invariant gradient. In \cite{Choe} or \cite{Zhu2005}, for $n=1$, the reader can find proof of boundedness and surjectivity of $P_{\alpha}$ from $L^{\infty}$ to $\mathcal{B}$. In \cite{Perala2012} and \cite{Perala2013}, Peralla found the exact value of the norm $\|P\|$ in $\mathbb{D}$ , while \cite{KalajMarkovic2014} contains a generalization of this result to $P_{\alpha}$ and $\mathbb{B}^n.$
  
  In \cite{KalajMarkovic2014}, the authors also have settled the problem of finding the exact value of $\|P_{\alpha}\|_{L^{\infty}(\mathbb{B}^n) \to \mathcal{B}(\mathbb{B}^n)}$ with a different norm on $\mathcal{B}.$ They obtained the two-sided estimate and conjectured
  that the norm is equal to the estimate from above. Using a new technique, we will obtain the appropriate series expansion of certain elliptic integral considered in \cite{KalajMarkovic2014} and the exact norm as the maximum of that series. We hope that technique can be used for a variety of  similar extremal problems. 
 
 Let us, first, recall that the Bloch space consists of functions $f$ analytic in $\mathbb{B}^n$ for which the following semi-norm is finite:
 $$ \|f\|_{\beta}:= \sup_{|z|<1} (1-|z|^2)|\nabla f(z)|,$$
 where $$\nabla f(z)=\Big(\frac{\partial f}{\partial z_1}(z),\dots,\frac{\partial f}{\partial z_n}(z)\Big).$$
 But, we can also define the semi-norm invariant with respect to the group $Aut(\mathbb{B}^n).$ For analytic $f$, the invariant gradient $\widetilde{\nabla} f(z)$ is defined by:
 $$\widetilde{\nabla} f(z)= \nabla (f\circ \varphi_z)(0),$$
 where $\varphi_z$ is an automorphism of the unit ball for which $\varphi_z(0)=z.$ We have
 $$ |\widetilde{\nabla}(f \circ \varphi)|=|(\widetilde{\nabla}f) \circ \varphi)|, $$
 exactly what we want. 
 Then the Bloch space can be described also as the space of all holomorphic functions $f$ for which 
 $$\|f\|_{\widetilde{\beta}}:=\sup_{|z|<1} |\widetilde{\nabla} f(z)|< \infty .$$
 
 Now, we can equip $\mathcal{B}$ with the norm $\|f\|_{\widetilde{\mathcal{B}}}:=|f(0)|+\|f\|_{\widetilde{\beta}}.$
 
 \subsection{Statement of the problem}
 
 In order to formulate the problem and the known result, we define the following function of 
 one real variable $t \in [0,\frac{\pi}{2}]:$
   
   $$  l(t)=(n+\alpha+1) \int_{\Balln} \frac{|(1-w_1)\cos t + w_2\sin t|}{|w_1-1|^{n+\alpha+1}} d v_{\alpha}.$$
   
   Kalaj and Marković  in \cite{KalajMarkovic2014} proved:
   
   \begin{theorem} For $\alpha>-1,$ $n>1,$ we have 
   $$ l(\frac{\pi}{2})= \frac{\pi}{2} l(0)= \frac{\pi}{2}C_{\alpha},  $$
   where $C_{\alpha}=\frac{\Gamma(n+\alpha+2)}{\Gamma^2(\frac{n+\alpha}{2}+1)}.$
   
   For the $\widetilde{\beta}-$semi-norm of the Bergman projection $P_{\alpha}$ we have:
   $$ \widetilde{C}_{\alpha}:=\|P_{\alpha}\|_{\widetilde{\beta}}= \max_{0 \leq t \leq \frac{\pi}{2}} l(t),   $$
   and 
   $$   \frac{\pi}{2}C_{\alpha} \leq \|P_{\alpha}\| \leq \frac{\sqrt{\pi^2 +4}}{2}C_{\alpha} . $$
   
   \end{theorem}
    They also conjectured that 
    $$ \|P_{\alpha}\|_{\widetilde{\beta}}= \frac{\pi}{2} C_{\alpha}. $$
    
    From these facts we can conclude that it is enough to prove that $l(t)$ attains its maximum in $t=\frac{\pi}{2}.$
    
    We will prove that this conjecture is true. This is contained in the following theorem.
    
    \begin{theorem} For $\alpha> -1$ and $n \geq 2$, we have  
    	$$ \|P_{\alpha}\|_{\widetilde{\beta}}= \frac{\pi}{2}\frac{\Gamma(n+\alpha+2)}{ \Gamma^2(\frac{n+\alpha}{2}+1)}
    	; \quad  \|P_{\alpha}\|_{\widetilde{\mathcal{B}}}= 1+ \frac{\pi}{2}\frac{ \Gamma(n+\alpha+2)}{ \Gamma^2(\frac{n+\alpha}{2}+1)}.$$
    	
    	Moreover, we have the following series expansion for function $l:$
    	
    	 $$l(t)= \frac{\Gamma(n+\alpha+2)}{\Gamma^2(\frac{n+\alpha}{2}+1)} \sin^2t\sum_{m=0}^{+\infty} \frac{\Gamaf{m+\frac{1}{2}} \Gamaf{m+\frac{3}{2}}}{\paren{m!}^2}\cos^{2m}t, \quad 0 < t \leq \frac{\pi}{2}$$
    	 
    	 and $l(t)$ is increasing in $ t \in [0,\frac{\pi}{2}].$

    	 \end{theorem}
     
     In the next section we give some preliminary facts which we need for the proof.
    
    \subsection{Hypergeometric functions}
    
    Here we recall some properties of hypergeometric functions.They are defined by
    $$ _2 F_1(a,b;c;z)= \sum_{k=0}^{+\infty} \frac{(a)_k (b)_k}{(c)_k} \frac{z^k}{k!}.     $$
    
    It converges for all $|z|<1 $, and, for $Re(c-a-b)>0$ also for $z=1.$ Here $(a)_k$ stands for $a(a+1) \dots (a+k-1), $ and $a$ is not negative integer.
    
    We will use the next theorem due to Gauss:
    $$ _2F_1(a,b;c;1)=\frac{\Gamma(c)\Gamma(c-a-b)}{\Gamma(c-a)\Gamma(c-b)}.     $$
    
   \section{Proof of the Theorem 2}

  Let us recall the integral representation of the constant $C_{\alpha}$. We start from the expression
\begin{align*}
   L(\xi_t):=& c_{\alpha} \int_{\Balln} \frac{\abs{\innerp{w-e_1}{\xi_t}} \paren{1 - \abs{w}^2}^{\alpha}}{\abs{1 - \innerp{w}{e_1}}^{n+\alpha+1}} \diff v_n(w)\\
   =& c_{\alpha} \int_{\Balln} \frac{\abs{(1-w_1) \cos t + w_2 \sin t}}{\abs{1-w_1}^{n+\alpha+1}} \paren{1-\abs{w}^2}^{\alpha} \diff v_n(w)
\end{align*}
where $\xi_t = e_1 \cos t + e_2 \sin t$,  $t \in \ccint{0}{\frac{\pi}{2}}$ and $c_{\alpha}$
is given in (1).

Let us fix $t \in [0,\frac{\pi}{2}]$. Changing coordinates with $A_tw=z$, where $A_t$ is a real  $n \times n$ orthogonal matrix
$$
\begin{pmatrix}
	\cos t & \sin t & 0 & \cdots & 0 \\
	-\sin t & \cos t & 0 & \cdots & \vdots \\
	0 & 0 & 1 & \cdots & 0 \\
	\vdots & \vdots & \vdots & \ddots & 0 \\
	0 & 0 & 0 & \cdots & 1
	
\end{pmatrix}$$

such that $A_t\xi_t = e_1$, we obtain
\begin{align*}
  L(\xi_t) =& c_{\alpha} \int_{\Balln} \frac{\abs{\innerp{A_tw-A_te_1}{e_1}} \paren{1-\abs{w}^2}^\alpha}{\abs{1-\innerp{A_t w}{A_t e_1}}^{n+\alpha+1}} \diff v_n(w)\\
 =& c_{\alpha} \int_{\Balln} \frac{\abs{\innerp{z-A_te_1}{e_1}} \paren{1-\abs{z}^2}^\alpha}{\abs{1-\innerp{z}{A_t e_1}}^{n+\alpha+1}} \diff v_n(z).
\end{align*}

Since $A_t e_1=\paren[\big]{\cos t, -\sin t, 0, \dots, 0}$, we have:
\begin{align*}
 L(\xi_t) =& c_{\alpha} \int_{\Balln} \frac{\abs{z_1 -
    \cos t} \paren{1-\abs{z}^2}^\alpha}{\abs{1-z_1 \cos t + z_2 \sin t}^{n+\alpha+1}} \diff v_n(z)\\ =&
     c_{\alpha} \int_{\Balln} \frac{\abs{z_1 -
    		\cos t} \paren{1-\abs{z}^2}^\alpha}{\abs{1-z_1 \cos t - z_2 \sin t}^{n+\alpha+1}} \diff v_n(z).
\end{align*}
  Now, as in \cite{Melentijevic2017}, we use Fubini's theorem:
\begin{align*}
\MoveEqLeft  L(\xi_t) = c_{\alpha} \int_{\Balln} \frac{\abs{z_1 - \cos t} \paren{1-\abs{z_1}^2 - \abs{z_2}^2 - \abs{z'}^2}^\alpha}{\abs{1-z_1 \cos t - z_2 \sin t}^{n+\alpha+1}} \diff v_n(z) \\
  &= c_{\alpha}\int\limits_{\Ball{2}} \frac{\abs{z_1 - \cos t} \diff v_2(z_1, z_2)}{\abs{1 - z_1 \cos t - z_2 \sin t}^{n+\alpha+1}} \int\limits_{\sqrt{1 - \abs{z_1}^2 - \abs{z_2}^2} \Ball{n-2}} \paren{1 - \abs{z_1}^2 - \abs{z_2}^2 - \abs{z'}^2}^{\alpha} \diff v_{n-2}(z');
\end{align*}
here $z=\paren{z_1, z_2, z'}$, $z' \in \setC^{n-2}$.

We make a substitution $z'=\lambda w,$ $\lambda=\sqrt{1 - \abs{z_1}^2 - \abs{z_2}^2}$, in the inner integral, which gives:
\begin{align*}
  \int_{\lambda \Ball{n-2}} \paren{\lambda^2 - \abs{z'}^2}^{\alpha} \diff v_{n-2}(z')
  = \lambda^{2 \alpha + 2n - 4} \int_{\Ball{n-2}} \paren{1 - \abs{w}^2}^{\alpha} \diff v_{n-2}(w).
\end{align*}

We easily find $\int_{\Ball{n-2}} \paren{1 - \abs{w}^2}^\alpha \diff v_{n-2}(w) =
k_{\alpha}=\pi^{n-2}\frac{\Gamma(\alpha+1)}{\Gamma(\alpha+n-1)}$, so
\[
L(\xi_t)=c_{\alpha} k_{\alpha} I(\cos t, \sin t), \]
 where
  \[ I(\cos t, \sin t)= \int_{\Ball{2}} \frac{\abs{z_1 - \cos t} \paren{1 -
    \abs{z_1}^2 - \abs{z_2}^2}^{n+\alpha-2}}{\abs{1 - z_1\cos t - z_2\sin t}^{n+\alpha+1}}
\diff v_2(z_1, z_2).
\]
  Now, the proof of Theorem 2 is reduced to proving monotonicity of $ I(\cos t, \sin t)$ as a function of $0 \leq t \leq \frac{\pi}{2}.$
  
Again, Fubini's theorem gives us:
\[
 I(\cos t, \sin t) =  \int\limits_{\mathbb{D}} \abs{z_1 - \cos t} \diff v(z_1)
  \int\limits_{\sqrt{1-\abs{z_1}^2} \mathbb{D}} \frac{\paren{1 - \abs{z_1}^2 -
      \abs{z_2}^2}^{n+\alpha-2}}{\abs{1 - z_1\cos t - z_2\sin t}^{n+\alpha+1}} \diff v(z_2).
\]

Next, we make substitution $z_2 = \sqrt{1 - \abs{z_1}^2} \rho e^{i \theta}$, $0 \leq \rho < 1$, $\theta \in
\ccint{0}{2 \pi}$:
\begin{align*}
\MoveEqLeft \int\limits_{\sqrt{1-\abs{z_1}^2} \mathbb{D}} \frac{\paren{1 - \abs{z_1}^2 - \abs{z_2}^2}^{n+\alpha-2}}{\abs{1 - z_1\cos t - z_2\sin t}^{n+\alpha+1}} \diff v(z_2) \\
&= \int_0^1 \diff \rho \int_0^{2\pi} \frac{\paren{1-\abs{z_1}^2}^{n+\alpha-2} \paren{1 - \rho^2}^{n+\alpha-2} \rho \paren{1-\abs{z_1}^2}}{\abs{1 - z_1\cos t - \sqrt{1-\abs{z_1}^2}\rho\sin t e^{i\theta}}^{n+\alpha+1}} \diff \theta \\
  &= \paren{1 - \abs{z_1}^2}^{n+\alpha-1} \int_0^1 \rho \paren{1 - \rho^2}^{n+\alpha-2} \diff \rho \int_0^{2\pi} \frac{\diff \theta}{\abs{1 - z_1\cos t - \sqrt{1-\abs{z_1}^2}\rho\sin t e^{i\theta}}^{n+\alpha+1}}. 
\end{align*}
  Next, we use Parseval's identity and Taylor's expansion of $(1-z)^{-\frac{n+\alpha+1}{2}}:$

\begin{align*}
   \MoveEqLeft \int_0^{2\pi} \frac{\diff \theta}{\abs{1 - z_1\cos t - \sqrt{1-\abs{z_1}^2}\rho\sin t e^{i\theta}}^{n+\alpha+1}} \\
  &= \frac{1}{\abs{1-z_1\cos t}^{n+\alpha+1}} \int_0^{2\pi} \frac{\diff \theta}{\abs*{1 - \frac{\sqrt{1-\abs{z_1}^2}\rho\sin t}{1-z_1\cos t} e^{i\theta}}^{n+\alpha+1}} \\
  &= \frac{2\pi}{\abs{1-z_1\cos t}^{n+\alpha+1}} \sum_{k=0}^{+\infty} \binom{\frac{n+\alpha+1}{2} + k - 1}{k}^{\mkern -5mu 2} \frac{\paren{1-\abs{z_1}^2}^k \rho^{2k} \sin^{2k}t}{\abs{1-z_1\cos t}^{2k}}.
\end{align*}

Therefore, we have:
\begin{align*}
  \MoveEqLeft \int\limits_{\sqrt{1-\abs{z_1}^2} \mathbb{D}} \frac{\paren{1 - \abs{z_1}^2 - \abs{z_2}^2}^{n+\alpha-2}}{\abs{1 - z_1\cos t - z_2\sin t}^{n+\alpha+1}} \diff v(z_2) \\
  &= 2\pi \sum_{k=0}^{+\infty} {\binom{\frac{n+\alpha+1}{2} + k - 1}{k}}^{\mkern -5mu 2} \frac{\paren{1-\abs{z_1}^2}^{k+n+\alpha-1} \sin^{2k}t}{\abs{1-z_1\cos t}^{2k+n+\alpha+1}} \int_0^1 \rho^{2k+1} \paren{1-\rho^2}^{n+\alpha-2} \diff \rho \\
  &= \pi \sum_{k=0}^{+\infty} {\binom{\frac{n+\alpha+1}{2} + k - 1}{k}}^{\mkern -5mu 2} \frac{\paren{1-\abs{z_1}^2}^{k+n+\alpha-1} \sin^{2k}t}{\abs{1-z_1\cos t}^{2k+n+\alpha+1}} \Betaf{k+1}{n+\alpha-1},
\end{align*}
and hence:
\[
I(\cos t, \sin t) = \pi \sum_{k=0}^{+\infty} \Betaf{k+1}{n+\alpha-1} {\binom{\frac{n+\alpha+1}{2} + k -
    1}{k}}^{\mkern -5mu 2} \sin^{2k}t \int_{\mathbb{D}} \frac{\abs{z_1 - \cos t}
  \paren{1-\abs{z_1}^2}^{k+n+\alpha-1}}{\abs{1-z_1\cos t}^{2k+n+\alpha+1}} \diff v(z_1).
\]

We calculate these integrals by changing coordinates with $z_1 = \frac{\cos t -
  \zeta}{1-  \zeta \conj{\cos t}} = \frac{\cos t - \zeta}{1- \zeta \cos t}$ (since $\cos t \in \setR$). Here, we assume $t>0.$ Then, we have:
\[ \zeta = \frac{\cos t - z_1}{1- z_1 \cos t}, \quad J_{\setR} =
  \frac{\paren{1-\cos^2t}^2}{\abs{1 -  \zeta \cos t}^4}. \]

Also, we need the following identities
\[ 1 - z_1\cos t = 1 - \frac{\cos  - \zeta}{1 -  \zeta \cos t} \cos t = \frac{1 -
    \cos^2t}{1- \zeta \cos t} \]
and 
\[ 1-\abs{z_1}^2 = \frac{\paren{1-\cos^2 t}
    \paren{1-\abs{\zeta}^2}}{\abs{1-\zeta \cos t}^2}. \]
Using the above substitution, we get:
\begin{align*}
  \MoveEqLeft \int_{\mathbb{D}} \frac{\abs{z_1 - \cos t}
  \paren{1-\abs{z_1}^2}^{k+n+\alpha-1}}{\abs{1-z_1\cos t}^{2k+n+\alpha+1}} \diff v(z_1) \\
  &= \int_{\mathbb{D}} \abs{\zeta} \frac{\paren{1-\cos^2t}^{k+n+\alpha-1} \paren{1-\abs{\zeta}^2}^{k+n+\alpha-1}}{\abs{1-\zeta \cos t}^{2k+2n+2\alpha-2}} \frac{\abs{1-\zeta \cos t}^{2k+n+\alpha}}{\paren{1-\cos^2t}^{2k+n+\alpha}} \frac{\paren{1-\cos^2t}^2}{\abs{1-\zeta \cos t}^4} \diff v(\zeta) \\
  &= \paren{1-\cos^2 t}^{1-k} \int_{\mathbb{D}} \abs{\zeta} \frac{\paren{1 - \abs{\zeta}^2}^{k+n+\alpha-1}}{\abs{1-\zeta \cos t}^{n+\alpha+2}} \diff v(\zeta).
\end{align*}

Passing to polar coordinates, we have:
\[
  \int_{\mathbb{D}} \abs{\zeta} \frac{\paren{1 - \abs{\zeta}^2}^{k+n+\alpha-1}}{\abs{1-\zeta \cos t}^{n+\alpha+2}} \diff v(\zeta) =
  \int_0^1 r^2 \paren{1-r^2}^{k+n+\alpha-1} \diff r \int_0^{2\pi} \frac{\diff
    \varphi}{\abs{1-r \cos t e^{i\varphi}}^{n+\alpha+2}},
\]
and then,again,  by Parseval's identity:
\[
  \int_0^{2\pi} \frac{\diff \varphi}{\abs{1-r e^{i\varphi}\cos t}^{n+\alpha+2}} =
  2\pi \sum_{m=0}^{+\infty} \binom{\frac{n+\alpha+2}{2} + m - 1}{m}^{\mkern -5mu 2} r^{2m} \cos^{2m}t ,
\]
thus:
\begin{align*}
  \MoveEqLeft \int_{\mathbb{D}} \abs{\zeta} \frac{\paren{1 - \abs{\zeta}^2}^{k+n+\alpha-1}}{\abs{1-\zeta \cos  t}^{n+\alpha+2}} \diff v(\zeta) \\
  &= 2\pi \sum_{m=0}^{\infty} \binom{\frac{n+\alpha+2}{2} + m - 1}{m}^{\mkern -5mu 2} \cos^{2m}t \int_0^1 r^{2m+2} \paren{1-r^2}^{k+n+\alpha-1} \diff r \\
  &= \pi \sum_{m=0}^{\infty} \binom{\frac{n+\alpha+2}{2} + m - 1}{m}^{\mkern -5mu 2} \cos^{2m}t \Betaf{m+\tfrac{3}{2}}{k+n+\alpha}.
\end{align*}
This gives:
\begin{multline*}
  I(\cos t, \sin t) = \pi^2 \sum_{k=0}^{+\infty} {\binom{\frac{n+\alpha+1}{2} + k - 1}{k}}^{\mkern
    -5mu 2} \Betaf{k+1}{n+\alpha-1} \sin^{2k}t \cdot \\
  \cdot \sum_{m=0}^{\infty} \binom{\frac{n+\alpha+2}{2} + m - 1}{m}^{\mkern -5mu 2}
  \Betaf{m+\tfrac{3}{2}}{k+n+\alpha} \cos^{2m}t \paren{1-\cos^2 t}^{1-k}.
\end{multline*}
Since $\paren{1-\cos^2t}^{1-k} \sin^{2k} t = 1-\cos^2 t$, we
conclude:
\begin{multline}
  I(\cos t, \sin t) = \pi^2 \paren{1-\cos^2t} \sum_{k,m=0}^{+\infty} {\binom{\frac{n+\alpha+1}{2} +
      k - 1}{k}}^{\mkern -5mu 2} \binom{\frac{n+\alpha+2}{2} + m - 1}{m}^{\mkern -5mu
    2} \cdot \\
 \cdot \Betaf{k+1}{n+\alpha-1} \Betaf{m+\tfrac{3}{2}}{k+n+\alpha} \cos^{2m} t.
\end{multline}

Now, we consider the function $\phi$ defined as
\[
  \phi(x) = \paren{1-x} \sum_{k,m=0}^{+\infty} {\binom{\frac{n+\alpha+1}{2} + k - 1}{k}}^{\mkern
    -5mu 2} \binom{\frac{n+\alpha+2}{2} + m - 1}{m}^{\mkern -5mu 2}
  \Betaf{k+1}{n+\alpha-1} \Betaf{m+\tfrac{3}{2}}{k+n+\alpha} x^m,
\]
for $0\leq x <1$ (because of condition $0 \leq \cos t < 1$!)

Using $\Betaf{\alpha}{\beta} = \frac{\Gamaf{\alpha}\Gamaf{\beta}}{\Gamaf{\alpha+\beta}}$ we have:
\[
  \phi(x) = \Gamaf{n+\alpha-1} \paren{1-x} \sum_{k,m=0}^{+\infty} {\binom{\frac{n+\alpha+1}{2} + k -
      1}{k}}^{\mkern -5mu 2} \binom{\frac{n+\alpha+2}{2} + m - 1}{m}^{\mkern -5mu 2}
  \frac{k!\, \Gamaf{m+\tfrac{3}{2}}}{\Gamaf{k+n+\alpha+m+\tfrac{3}{2}}} x^m.
\]

Let us sum, over $k$, the terms which depend on $k$:
\begin{align*}
  \MoveEqLeft \sum_{k=0}^{+\infty} {\binom{\frac{n+\alpha+1}{2} + k - 1}{k}}^{\mkern -5mu 2} \frac{k!}{\Gamaf{k+n+\alpha+m+\tfrac{3}{2}}} \\
  &= \sum_{k=0}^{+\infty} \frac{\paren[\big]{\frac{n+\alpha+1}{2}+k-1}^2 \cdot \dotsm \cdot \paren[\big]{\frac{n+\alpha+1}{2}+1}^2 \paren[\big]{\frac{n+\alpha+1}{2}}^2}{k!\, \Gamaf{k+m+n+\alpha+\tfrac{3}{2}}} \\
  &= \sum_{k=0}^{+\infty} \frac{\paren[\big]{\frac{n+\alpha+1}{2}}_k \paren[\big]{\frac{n+\alpha+1}{2}}_k}{k!\, \paren[\big]{n+k+\alpha+m+\tfrac{1}{2}} \cdot \dotsm \cdot  \paren[\big]{n+\alpha+m+\tfrac{3}{2}} \Gamaf{n+\alpha+m+\tfrac{3}{2}}} \\
  &= \frac{1}{\Gamaf{n+\alpha+m+\tfrac{3}{2}}} \sum_{k=0}^{+\infty} \frac{1}{k!} \frac{\paren[\big]{\frac{n+\alpha+1}{2}}_k \paren[\big]{\frac{n+\alpha+1}{2}}_k} {\paren[\big]{n+\alpha+m+\tfrac{3}{2}}_k}.
\end{align*}
We recognize that the last sum is $\pFq{2}{1}{\frac{n+\alpha+1}{2}, \frac{n+\alpha+1}{2};
  n + \alpha + m + \frac{3}{2}; 1}$, and by Gauss's theorem this is equal to
\[ \frac{\Gamaf{n+\alpha+m+\frac{3}{2}} \Gamaf{m+\frac{1}{2}}}{\Gamma^2
    \paren{m+1+\frac{n+\alpha}{2}}}. \]
Hence, the double sum in (2) is equal to:
\[
  \sum_{m=0}^{+\infty} \binom{\frac{n+\alpha+2}{2} + m - 1}{m}^{\mkern -5mu 2}
  \frac{\Gamaf{m+\frac{3}{2}} \Gamaf{m+\frac{1}{2}}}{\Gamma^2
    \paren{m+1+\frac{n+\alpha}{2}}} x^m.
\]

Note that
\[
\binom{\frac{n+\alpha+2}{2} + m - 1}{m}^{\mkern -5mu 2} = \frac{1}{\paren{m!}^2}
\paren[\big]{\tfrac{n+\alpha+2}{2} + m - 1}^2 \cdot \dotsm \cdot
\paren[\big]{\tfrac{n+\alpha+2}{2}}^2 = \frac{1}{\paren{m!}^2} \frac{\Gamma^2
  \paren{\frac{n+\alpha+2}{2} + m}}{\Gamma^2 \paren{\frac{n+\alpha+2}{2}}},
\]
and hence
\[
\phi(x) = \frac{\Gamaf{n+\alpha-1}}{\Gamma^2 \paren{\frac{n+\alpha}{2} + 1}} \paren{1-x}
\sum_{m=0}^{+\infty} \frac{\Gamaf{m+\frac{1}{2}} \Gamaf{m+\frac{3}{2}}}{\paren{m!}^2}
x^m, \quad 0 \leq x < 1.
\]

Let us denote
\[
a_m = \frac{\Gamaf{m+\frac{1}{2}} \Gamaf{m+\frac{3}{2}}}{\paren{m!}^2}.
\]
It is easily verified that $a_m$ is strictly decreasing in $m \geq 0:$
\[
\frac{a_{m+1}}{a_m} = \frac{\Gamaf{m + \frac{3}{2}} \Gamaf{m+\frac{5}{2}}
  \paren{m!}^2}{\Gamaf{m + \frac{1}{2}} \Gamaf{m+\frac{3}{2}}
  \paren[\big]{\paren{m+1}!}^2} = \frac{\paren{m+\frac{1}{2}}
  \paren{m+\frac{3}{2}}}{\paren{m+1}^2} < 1.
\]
In particular,  $a_m \leq a_0$.

Then, we may conclude
\[ (1-x) \sum_{m=0}^{+\infty} a_m x^m \leq (1-x) a_0 \sum_{m=0}^{+\infty} x^m = a_0,\]
that is  $\phi(x) \leq \phi(0)$.
Moreover, $\phi(x)$ is decreasing, since we can write it in the following form:
$$  \phi(x)= \frac{\Gamaf{n+\alpha+1}}{\Gamma^2({\frac{n+\alpha}{2}+1})} \big(a_0 + \sum_{m=1}^{+\infty} (a_m -a_{m-1}) x^{m} \big). $$

This is the crux of the proof of our Theorem. 

So, the best constant is  $\widetilde{C}_{\alpha}=(1+n+\alpha) c_{\alpha} k_{\alpha} \pi^2 \frac{\Gamma(\alpha+n-1)}{\Gamma^2(\frac{n+\alpha}{2}+1)} \Gamma(\frac{1}{2})\Gamma(\frac{3}{2})$, i.e. 
 $$\widetilde{C}_{\alpha}= \frac{\pi}{2} \frac{\Gamma(\alpha+n+2)}{\Gamma^2(\frac{n+\alpha}{2}+1)}.$$
 
 According to \cite{KalajMarkovic2014}, for  $\xi=(1,0,0 \dots ,0)$ we have $l(0)=\frac{2}{\pi} \widetilde{C}_{\alpha}$. This can be also obtained from the above series by letting $x$  tends to $1.$

  Finally, all these computations give us 
 $$l(t)= (1+n+\alpha) c_{\alpha} k_{\alpha} \pi^2 \frac{\Gamma(\alpha+n-1)}{\Gamma^2(\frac{\alpha+n}{2}+1)} \sin^2t\sum_{m=0}^{+\infty} \frac{\Gamaf{m+\frac{1}{2}} \Gamaf{m+\frac{3}{2}}}{\paren{m!}^2}
 \cos^{2m} t=$$
 $$ \frac{\Gamma(n+\alpha+2)}{\Gamma^2(\frac{n+\alpha}{2}+1)} \sin^2t\sum_{m=0}^{+\infty} \frac{\Gamaf{m+\frac{1}{2}} \Gamaf{m+\frac{3}{2}}}{\paren{m!}^2}\cos^{2m}t, \quad 0 < t \leq \frac{\pi}{2},$$
 and $l(t)$ is increasing in $t \in [0,\frac{\pi}{2}].$ (Because $\phi(x)$ is decreasing and $l(t)=\phi(\cos^2t).$)

 This concludes the proof of Theorem 2. 
 
 \textbf{Acknowledgements}. The author is partially supported by MPNTR grant, no. 174017, Serbia.

\end{document}